\documentclass{amsart}
\usepackage{latexsym}

\title{Harnack Inequality for a class of degenerate elliptic operators}
\author{Jos\'e D. Fernandes \and Jorge Groisman 
\and Severino T. Melo
}
\date{}
\newtheorem{thm}{Theorem}
\newtheorem{prop}{Proposition}
\newtheorem{lemma}[prop]{Lemma}
\newtheorem{obs}[prop]{Remark}
\newtheorem{defn}[prop]{Definition}
\begin{document}

\newcommand{\rn}{{\mathbb R}^{n}}
\newcommand{\rrm}{{\mathbb R}^{m}}
\newcommand{\rN}{{\mathbb R}^{N}}
\newcommand{\op}{operator}
\newcommand{\ops}{operators}
\newcommand{\pb}{\overline{p}}
\newcommand{\cqd}{\hfill$\Box$}
\newcommand{\dx}{d\!x}
\newcommand{\dz}{d\!z}
\newcommand{\nl}{\nabla_{\lambda}}
\newcommand{\sumij}{\sum_{i,j=1}^{N}}
\newcommand{\ho}{H(\Omega)}
\newcommand{\hoo}{H_{\circ}(\Omega)}
\newcommand{\hq}{H(Q)}
\newcommand{\hoq}{H_{\circ}(Q)}
\newcommand{\ap}{{\em (}}
\newcommand{\fp}{{\em )}}
\newcommand{\prf}{{\em Proof}: }
\newcommand{\hmb}{H_{_{M,d}}}
\newcommand{\fkj}{f_{k_{j}}}
\newcommand{\lip}{\mbox{{\tt Lip}}}
\newcommand{\lipo}{\mbox{{\tt Lip}}_{\circ}}

\begin{abstract}
We prove a Harnack inequality for a class of two-weight degenerate elliptic operators. The metric distance is induced by continuous Grushin-type vector fields. It is not know whether there exist cutoffs fitting the metric balls. This obstacle is bypassed by means of a covering argument that allows the use of rectangles in the Moser iteration. 
\end{abstract}
\maketitle

\section{Introduction}

Perhaps inspired by David and Semmes' work \cite{DS}, Franchi, Gutierrez and 
Wheeden proved in \cite{FGW1}\ a very deep generalization of the classical
Sobolev-Poincar\'e inequality, unifying several other previous 
results. The importance of Sobolev-Poin\-ca\-r\'e-type inequalities to 
the study of elliptic equations has been well known for decades \cite{GT}. 
In particular, the so-called {\em Moser iteration}\
technique \cite{M0,Mp5,M1}\ still is the basis upon which are built 
more recent proofs of Harnack-type inequalities for non-negative 
solutions of degenerate elliptic equations \cite{BM,ChW,italianas,FKS,FL,FSSC,G}. 

The main result in \cite{FGW1}\ thus paved the way for the proof of 
a more general Harnack inequality. Indeed, in \cite{FGW}, Theorem II, the same authors 
stated a result which has as particular cases the Harnack inequalities 
proven in \cite{ChW}\ and \cite{FKS}. As they pointed out, that new version
would apply to solutions of the equation
\begin{equation}
\label{fgwex}
\frac{\partial}{\partial x}\left[
(|x|^{\sigma +1}+|y|)^{\frac{\kappa}{\sigma +1}}
\frac{\partial f}{\partial x}\right] 
+
\frac{\partial}{\partial y}\left[
(|x|^{\sigma +1}+|y|)^{\frac{\kappa}{\sigma +1}}
|x|^{\sigma}\frac{\partial f}{\partial y}\right] 
=0
\end{equation}
in an open set $\Omega\subset{\mathbb R}^2$\ containing the origin, 
with $\kappa$\ and $\sigma$\ arbitrary positive numbers. None of 
the other available results includes this example.

The proof of Theorem II in \cite{FGW}, however, is not complete. It depends
on the (not proven) existence of certain cut-off functions fitting the metric balls defined by the operator. It is easy to construct (see our 
Proposition~\ref{test}, below) cutoffs which are identical to one or nonzero 
not on metric balls, but on certain ``rectangles'' which are products of Euclidean balls with variable ratio of the radii. 
If one insists in using balls contained or containing those rectangles, there 
remains a gap between the two balls which provokes an explosion
of the constants that appear in the iteration process.

In this paper, we prove Theorem II of \cite{FGW}\ without using cutoffs 
addapted to balls, applying instead a covering 
technique, based on a theorem in \cite{CoW}, already used in the study 
of degenerate parabolic equations by the first author \cite{JF}. The 
building block of the Moser iteration used here turns out to be not exactly a
Sobolev-Poincar\'e inequality, but rather its consequence stated in 
Theorem~\ref{D}; which is a Sobolev-Poincar\'e inequality for rectangles, 
with the one on the right $\epsilon$\ times larger than the one on 
the left and with a negative power of $\epsilon$\ on the right. The main point
of Section~\ref{moserze}\ is to show that a sequence $\epsilon_k$\ can be 
chosen in such a way that the iteration converges. We show that the 
Moser-type iteration designed by Chanillo and Wheeden in \cite{ChW}\ also 
works in this context. Propositions which are straighforward addaptions of results in \cite{ChW}\ are stated here without proof.  

We will assume as a hypothesis that the 
Sobolev-Poincar\'e inequality we need is true, without explicitly stating 
Franchi, Gutierrez and Wheeden's Theorem~I of \cite{FGW1}, which is 
nonetheless our main motivation (since it provides the main example). 
One important aspect of that theorem is that 
it allows the presence of two (possibly non-comparable and non-Muckenhoupt) weights in the ellipticity condition. 

The existence of cutoffs suitable to the study of regularity properties of weak solutions of degenerate elliptic equations has been independently proven by Franchi, Serapioni and Serra Cassano \cite{FSSC}, and by Garofalo and Nhieu \cite{GN}. Their results would apply in our context, however, 
only if we required that the function $\lambda$, defined in our Section~\ref{primsec}, be Lipschitz continuous (for the  operator in (\ref{fgwex}),  the natural choice of $\lambda$\ would be 
$\lambda(x)=|x|^{\sigma}$, $\sigma>0$). Under this additional assumption, Theorem~1.3 in \cite{GN}, or Proposition~2.9
in \cite{FSSC}\ (together with, for example, the composition argument in the proof of Theorem~1.5 in \cite{GN}), would
imply the existence of the test functions needed for the proof of Theorem~II in \cite{FGW}\ to work. 

A different approach was taken by Biroli and Mosco \cite{BM}. Within a very general framework, they proved the 
existence of cutoffs which satisfy, instead of a pointwise estimate (as in \cite{GN}, Theorem~1.5, for example), a weaker 
requirement, in integral form (\cite{BM}, Proposition~3.3). That also suffices for the proof of Harnack-type inequalities (Theorem~1.1 in \cite{BM}; Theorem~1 in \cite{G}). Working directly with the bilinear form defined by the elliptic operator, they did not have to 
to deal with the  regularity of the vector fields usually used to define the metric. 

%
%
%
%
\section{Preliminaries and statement of the main result}
\label{primsec}

The \ops\ considered in this paper are of type
\begin{equation}
\label{L}
Lf\,=\,\sum_{i,j=1}^{N}\frac{\partial}{\partial z_i}\left(a_{ij}(z)
                \frac{\partial f}{\partial z_j}\right)\,,
\end{equation}
where $z=(z_1,\cdots,z_N)=(x_1,\cdots,x_n,y_1,\cdots,y_m)\in{\mathbb R}^N=
\rn\times{\mathbb R}^m,$\ 
the matrix $A=(\!(a_{ij})\!)$\ is symmetric and the functions $a_{ij}$\ are real, 
measurable and satisfy the (degenerate) ellipticity condition
\begin{equation}
\label{ell}
v(z)(|\xi|^2+\lambda(x)^2|\eta|^2)
\leq
\sumij a_{ij}(z)\zeta_{i}\zeta_{j}
\leq
u(z)(|\xi|^2+\lambda(x)^2|\eta|^2),
\end{equation}
for all $\zeta=(\xi,\eta)\in{\mathbb R}^n\times{\mathbb R}^m,$\ 
with the functions $\lambda,$\ $u$\ and $v$\ non-negative and satisfying
several hypotheses which are especified in what follows.

Throughout this paper, $aB$\ will denote, for $a>0$\ and $B$\ a ball in some 
metric space, another ball with the same center and $a$-times the radius as $B.$

We require that the function $\lambda,$\ defined on $\rn,$\ satisfy:

\begin{description}

\item[{\bf H1}] It is non-negative, continuous, and vanishes 
possibly only on a set of isolated points.

\item[{\bf H2}] It is {\em doubling}\ with respect to the Euclidean metric
and the Lebesgue measure, with doubling constant $C_1$; i.e., 
\[
\int_{2B_e}\lambda(x)\dx\leq C_1\int_{B_{e}}\lambda(x)\dx, 
\]
for every Euclidean ball $B_e\subset\rn.$

\item[{\bf H3}] There exists a constant $C_2$\ such that         \[
\sup_{x\in B_e}\lambda(x)\leq C_2\frac{1}{|B_e|}\int_{B_e}\lambda(x)\dx\,,
\]
for every Euclidean ball $B_e\subset\rn,$\ with $|\cdot|$\ denoting the 
Lebesgue measure.

\end{description}


\begin{defn}
\label{cubos}
Given $z_{\circ}=(x_{\circ},y_{\circ})\in{\mathbb R}^n\times{\mathbb R}^m=
{\mathbb R}^N$\ and $r>0,$\ we define
\[
\Lambda(z_{\circ},r)\,=\,\sup_{\{x;|x-x_{\circ}|<r\}}\lambda(x)
\]
and denote 
\[
Q(z_{\circ},r)\,=\,\{(x,y)\in{\mathbb R}^n\times{\mathbb R}^m;
\,|x-x_{\circ}|<r\,,\ |y-y_{\circ}|<r\Lambda(z_{\circ},r)\,\}.
\]
If $Q=Q(z_{\circ},r)$\ and $t>0,$\ $tQ$\ will denote $Q(z_{\circ},tr).$
\end{defn}

\begin{obs}
\label{2}
{\em
If follows from ({\bf H2}) and ({\bf H3}) that 
\begin{equation}
\label{smdob}
\Lambda(z_{\circ},2r)\leq
\frac{C_1C_2}{2^{n}}\Lambda(z_{\circ},r)
\end{equation}
for all $z_{\circ}\in\rN$\ and all $r>0;$\ and, hence, $C_1C_2\geq 2^n$\ must hold.
}
\end{obs}

\begin{lemma}
\label{psdmet}
If $z\in Q(z_{\circ},r)$\ and $w\in Q(z,s),$\ then $w\in Q(z_{\circ},r+s).$
\end{lemma}

\begin{defn} An absolutely continuous curve in ${\mathbb R}^N$\ is   
{\em subunit}\ if, for every $\zeta=(\xi,\eta)\in{\mathbb R}^N$\ 
and for almost every $t$\ in its domain, we have
\[
\langle\gamma^\prime(t),\zeta\rangle^2
\leq|\xi|^2+\lambda(\gamma(t))^2|\eta|^2\,,
\]
with $\langle\cdot,\cdot\rangle$\ denoting the usual inner-product 
of ${\mathbb R}^N.$\ Given $z$\ and $w$\ 
in ${\mathbb R}^N,$ let $\rho(z,w)$\ denote the infimum of all $T\geq 0$\
such that there is a subunit curve joining the two points
with domain $[0,T].$
\end{defn}

The function $\rho$\ corresponds to the metric on 
${\mathbb R}^N$\ associated to the {\em Grushin-type}\ vector fields 
$
\frac{\partial}{\partial x_{1}},\cdots,
\frac{\partial}{\partial x_{n}}, 
\lambda(x)\frac{\partial}{\partial y_{1}},\cdots,
\lambda(x)\frac{\partial}{\partial y_{m}}
$ 
in a way which has by now become standard \cite{FL1}.
If $\lambda$\ is smooth and does not vanish, one can see that 
$\rho$\ is equal to the geodesic distance associated to the Riemannian
metric 
$
d\!s^2=\sum_{i=1}^{n}\dx_i^2+\lambda(x)^{-2}\sum_{j=1}^{m}d\!y_j^2.
$

An elementary proof of the following proposition can be given. 
For a somewhat different but closely related result, we refer to \cite{F1}. 

\begin{prop}\label{doble}
The function $\rho$ above defines a metric on ${\mathbb R}^N$\ and
there exists a constant $b,$\ depending only on $n$\ and $m,$\ such that 
the double inclusion
\begin{equation}
\label{dbl}
Q(z_{\circ},r/b)\,\subseteq B(z_{\circ},r) \,\subseteq Q(z_{\circ},br)
\end{equation}
holds for every $z_{\circ}\in{\mathbb R}^N$\ and $r>0,$\
where $B(z_{\circ},r)$\ denotes the ball with respect to this new metric 
with center $z_{\circ}$\ and radius $r.$
\end{prop}

\begin{obs}
{\em Only ({\bf H1}) is required for the proof of Proposition~\ref{doble}. 
As shown in Proposition 2.1.1 of \cite{jorge}, one may take 
$b=\max\{3,\sqrt{m},\sqrt{n}\}.$\ 
Proposition~\ref{doble}\  and ({\bf H1}) imply that the metric 
$\rho$\ induces in $\rN$\ its usual topology.
}
\end{obs}

We require that $u$\ and $v$\ be {\em weights} on $\rN,$\ (non-negative 
non-trivial locally integrable functions), which are 
{\em doubling}\ with respect to the $\rho$-metric and the Lebesgue measure, 
i.e., such that there are positive constants $C_3$\ and $C_4,$\ with
\begin{equation}
\label{dob}
\int_{2B}u(z)\dz\leq C_3\int_{B}u(z)\dz\ \mbox{and}\  
\int_{2B}v(z)\dz\leq C_4\int_{B}v(z)\dz, 
\end{equation}
holding for all $\rho$-balls $B.$\ For every measurable $E\subseteq\rN,$\
we will denote by $u(E)$\ and $v(E)$\ the integrals over $E$\ of $u$\
and $v,$\ respectively. Notice that (\ref{dob}) and Proposition~\ref{doble}\ 
imply that $u(E)$\ and $v(E)$\ are positive if $E$\ has non-empty interior.

For every locally integrable function $g,$\
we will denote by $m_{E}(g)$\ \label{pagina}
the $u$-average $u(E)^{-1}\int_{E}gu.$

Last we state the strongest hypothesis we impose on 
$u,$\ $v$\ and $\lambda$: that the 
following Sobolev-Poincar\'e inequality  holds. For sufficient 
conditions for its validity see, for example, the papers 
\cite{ChW0,FKS,F1,FGW1,FL1,GNC,J,SW}\ and their references.

\begin{description}
\item[{\bf SP}]
There exist $q>2$\ and $C_5>0,$\ constants depending only on $u,$\ $v,$\
$\lambda,$\ $n$\ and $m,$\ such that the inequality
\[
\left[\frac{1}{u(B)}\int_{B}|g(z)-m_{B}(g)|^q u(z)\dz\right]^{\frac{1}{q}}
\leq C_5 r
\left[\frac{1}{v(B)}\int_{B}|\nl g(z)|^2 v(z)\dz\right]^{\frac{1}{2}}
\]
holds 
for every Lipschitz continuous function $g$\ and every ball $B$\ with respect
to the metric $\rho$\ induced by $\lambda,$\ with $r$\ denoting the radius of 
$B,$\ and $\nl g$\ denoting the vector field
\[
\nl g(z)=\left(
\frac{\partial g}{\partial x_{1}}(z),\cdots,
\frac{\partial g}{\partial x_{n}}(z), 
\lambda(x)\frac{\partial g}{\partial y_{1}}(z),\cdots,
\lambda(x)\frac{\partial g}{\partial y_{m}}(z)
\right).
\]
\end{description}

Weak solutions of $Lf=0$\ in a bounded open set $\Omega\subset\rN$\ 
are defined (as in \cite{FGW}) in $\ho,$\ the completion 
of the space $\lip(\overline{\Omega})$\ of the Lipschitz continuous 
functions on $\overline{\Omega},$\ 
the closure of $\Omega,$\ with respect to the norm
\begin{equation}
\label{hilbert}
||f||_{H}^2=
\sumij\int_{\Omega} a_{ij}(z)\frac{\partial f}{\partial z_{i}}(z)
\frac{\partial f}{\partial z_{j}}(z)\dz
+\int_{\Omega}f(z)^2u(z)\dz.
\end{equation}
Using (\ref{ell}) and (\ref{dob}),   
one can show, similarly as in \cite{ChW}, that the equation above indeed 
defines a norm. 
Moreover, if we denote by $\hoo$\ the closure in $\ho$\ of the space 
$\lipo(\Omega)$\ of the 
Lipschitz continuous functions of compact support in $\Omega,$\ it can be 
proven, and for that ({\bf SP}) is required, that the bilinear form
$a_{\circ}$\ on $\lipo(\Omega),$
\[
a_{\circ}(f,g)\,=\,
\sumij\int_{\Omega} a_{ij}(z)\frac{\partial f}{\partial z_{i}}(z)
\frac{\partial g}{\partial z_{j}}(z)\dz,
\]
induces on $\hoo$\ an inner-product whose corresponding norm is equivalent to
$||\cdot||_{H}.$

\begin{defn} {\em An element $f\in\ho$\ is a weak solution of $Lf=0$\ if
$a_{\circ}(f,\theta)=0$\ for all $\theta\in\hoo.$}
\end{defn}

Applying Lax-Milgram's Theorem, existence and uniqueness of a suitably 
defined weak version of the Dirichlet problem on $\Omega$\ can be proven,
in exactly the same way as in \cite{ChW}. 

We still need two more definitions. The inequality
$
\int_{\Omega}f(x)^2u(x)\dx\,\leq\,||f||_{H}^2
$\ 
follows from (\ref{ell}) and the definition of $||\cdot||_{H}.$\ A natural 
mapping $\ho\rightarrow L^2(\Omega,u(z)\dz),$\ 
$f\mapsto\tilde{f},$\ is then defined. We
stress we are not claiming that this is an injection, even though that could
be proven under additional hypotheses. Finally, we will call an 
$f\in\ho$\ non-negative, and denote this by $f\geq 0,$\ if there 
is a sequence of non-negative functions $f_k\in\lip(\overline{\Omega})$\ 
converging to $f$\ in $\ho.$\ 

\begin{obs}
\label{rstr}
{\em
If $U\subset\Omega$\ is open and $f\in\ho$\ is a weak solution of $Lf=0$\
in $\Omega,$\, the restriction $f\!\!\mid_{U}\in H(U)$\ is then a weak solution 
of $Lf=0$\ in $U.$\ We also have $\tilde{f}\!\!\mid_{U}=\tilde{f\!\!\mid_{U}}.$}
\end{obs}
We are ready to state our main result.
%
%
\begin{thm} 
\label{main}
Suppose that $\lambda$\ satisfies \ap{\bf H1}\fp, 
\ap{\bf H2}\fp\ and \ap{\bf H3}\fp, $u$\ and $v$\ are doubling 
weights and also that \ap{\bf SP}\fp\ holds. 
Then there is a constant $K,$\ depending only on 
 $C_1,$\ $C_2$\, $C_3,$\ $C_4,$\ $C_5,$\ $q,$\ $m$\ and $n,$\ 
such that, if $\Omega$\ is a bounded open subset of $\rN$\ and 
$f\in\ho$\ is a non-negative weak solution of $Lf=0,$\ with $L$\ satisfying 
\ap\!\ref{L}\fp\ and \ap\!\ref{ell}\fp, then
%
%
\begin{equation}
\label{harn}
\mbox{\tt ess\,sup}_B\,
\tilde{f}\,\leq\, e^{K\mu}\,
\mbox{\tt ess\,inf}_B\,\tilde{f},
\end{equation}
for every $\rho$-ball $B$\ such that 
$2b^{4}B\subseteq\Omega,$\ where $\mu=u(B)^{\frac{1}{2}}v(B)^{-\frac{1}{2}}.$
\end{thm}         

%
%
%
%
%
%
%
%

\section{Application of a covering technique} 
\label{seccov}

All hypotheses of Theorem~\ref{main} are assumed to be true for the rest of the 
paper, even if not explicitly. By a ``constant'' we will 
always mean a positive number which may depend only on the constants that
arise in the hypotheses of Theorem~\ref{main}: $C_1,$\ $C_2$\, $C_3,$\ 
$C_4,$\ $C_5,$\ $q,$\ $m$\ and $n.$\ We start with a Sobolev-Poincar\'e 
inequality  for the rectangles $Q$\ of Definition~\ref{cubos}. 

\begin{prop} 
\label{P1}
There exists a constant $C_6$\ such that 
\[
\left[\frac{1}{u(Q)}\int_{Q}|g(z)|^q u(z)\dz\right]^{\frac{1}{q}}
\leq 
C_6 r\left[\frac{1}{v(Q)}\int_{b^{2}Q}|\nl g(z)|^2 v(z)\dz\right]^{\frac{1}{2}}
\]
\begin{equation}
\label{spc}
+
\left[\frac{1}{u(Q)}\int_{b^2Q}g(z)^2u(z)\dz\right]^{\frac{1}{2}}
\end{equation}
holds for every Lipschitz function $g$\ and every $Q=Q(z,r),$\ 
where $q>2$\ is the constant provided by \ap{\bf SP}\fp.
\end{prop}

\prf 
Using (\ref{dbl}), we see that 
\[
\left[\frac{1}{u(Q)}\int_{Q}|g(z)-m_{bB}(g)|^q u(z)\dz\right]^{\frac{1}{q}}
\]
is bounded by
\[
\left[\frac{u(bB)}{u(\frac{1}{b}B)}\frac{1}{u(bB)}
\int_{bB}|g(z)-m_{bB}(g)|^q u(z)\dz\right]^{\frac{1}{q}}.
\]
Using that $u$\ is doubling and the inequality ({\bf SP}) for the ball $bB,$\ 
we get: 
\begin{equation}
\label{spc0}
\left[\frac{1}{u(Q)}\int_{Q}|g-m_{bB}(g)|^q u\right]^{\frac{1}{q}}
\leq C_6 
r\left[\frac{1}{v(Q)}\int_{b^{2}Q}|\nl g|^2 v\right]^{\frac{1}{2}},
\end{equation}
with $C_6=bC_3^{\frac{l}{q}}C_5,$\ where $l$\ is an integer
such that $b^2<2^l.$

To prove (\ref{spc}), we start by applying to $g=[g-m_{bB}(g)]+m_{bB}(g)$\ 
the triangle inequality in $L^p(Q,u(z)\dz),$\ followed by (\ref{spc0}), then 
by the Cauchy-Schwarz inequality for $L^2(Q,u(z)\dz)$\ and finally (\ref{dbl}).
\cqd

We will call a metric space {\em homogeneous}\ if it can be equipped with 
a Borel 
measure $\nu$\ such that $\nu(2B)\leq D\nu(B)$\ for every ball $B,$\ for some
{\em doubling-factor}\ $D.$\ The following proposition is a particular case 
of Theorem~1.2 of \cite{CoW}.

\begin{prop}
\label{cow}
If $\{B(x,r)\}$\ is a family of balls of constant radius covering a subset 
$E$\ of a homogeneous metric space X, then there is a finite sub-family 
$\{B(x_i,r);$\ $i=1,\cdots,m\}$\ of disjoint balls such that 
$\{B(x_i,4r);i=1,\cdots,m\}$\ still covers $E.$
\end{prop}

\begin{prop}
\label{homo}
The metric space $(\rN,\rho)$\ is homogeneous.
\end{prop}
\prf Let $z_{\circ}=(x_{\circ},y_{\circ})\in\rn\times\rrm$\ and $r>0$\ 
be given. By (\ref{smdob}), we have
\begin{equation}
\label{svdob}
\Lambda(z_{\circ},t)\,\leq\, C_7^l\Lambda(z_{\circ},\frac{t}{2^{l}})
\end{equation}
for every non-negative integer $l$\ and every $t>0,$\ with 
$C_7=2^{-n}C_{1}C_{2}.$\ Using Proposition~\ref{doble}, we then get
\[
|B(z_{\circ},2r)|\leq\omega_n\omega_m(2br)^N\Lambda(z_{\circ},2br)^m
\leq C_7^{ml}(2b^2)^N|Q(z_{\circ},r/b)|\,,
\]
if $l$\ is chosen so that $2b^2\leq 2^l,$\ with $\omega_k$\ denoting the 
volume of the unit ball in ${\mathbb R}^r.$\ Since 
$Q(z_{\circ},r/b)\subseteq B(z_{\circ},r),$\ this shows that 
the Lebesgue measure is doubling with doubling-factor $C_7^{ml}(2b^2)^N.$
\cqd

\begin{prop}
\label{cover}
Given $z\in\rN$\ and $0<r<s,$\ there exist $z_1,\cdots,z_p$\ in $Q(z,s),$\ 
such that the family $\{Q(z_{1},r),\cdots,Q(z_{p},r)\}$\ covers $Q(z,s),$\ 
with $Q(z_j,\frac{r}{4b^{2}})$\ and $Q(z_k,\frac{r}{4b^{2}})$\ disjoint when 
$j\neq k.$\ Moreover, there are constants $\beta$\ and $C_8$\ such that
\begin{equation}
\label{cota}
p\,\leq\,C_8\left(\frac{s}{r}\right)^{\beta}\,.
\end{equation}
\end{prop}
\prf The first statement of this proposition follows straightforwardly from
Proposition~\ref{doble}, Proposition~\ref{cow}\ (with $\frac{r}{4b}$\ 
replacing $r$) and Proposition~\ref{homo}. 
In order to prove (\ref{cota}), let us first remark that there is a constant 
$\beta$\ such that the inequality
\begin{equation}
\label{invdob}
|Q(w,\theta t)|\,\geq\,C_7^{-m}\theta^{\beta}|Q(w,t)|
\end{equation}
holds for all $0<\theta<1,$\ $t>0$\ and $w\in\rN.$\ Indeed, let $\beta$\ be 
defined by $\beta=N+m\log C_7/\log 2.$\ Using 
$|Q(w,t)|=\omega_n\omega_mt^N\Lambda(w,t)^m,$\ 
we get (\ref{invdob}) by applying (\ref{svdob})
to the integer $l$\ such that $\theta/2<2^{-l}\leq\theta.$\ It follows
from Remark~\ref{2}\ that $C_7\geq 1$\ and thus $\beta$\ is positive.

By Lemma~\ref{psdmet}, and since $s+\frac{r}{4b^2}<(b^2+1)s,$\ 
each $Q_j=Q(z_j,\frac{r}{4b^{2}})$\ is contained in 
$Q(z,(b^2+1)s).$\ Since the $Q_j$'s are mutually disjoint, we have:
\begin{equation}
\label{rro}
|Q(z,(b^2+1)s)|\,\geq\,\sum_{j=1}^{p}|Q(z_j,\frac{r}{4b^{2}})|\,.
\end{equation}
Now let us apply (\ref{invdob}) to $w=z_j,$\ $t=(2b^2+1)s$\ and 
$\theta=r/(8b^4s+4b^2s).$\ We get:
\begin{equation}
\label{boi}
|Q(z_j,\frac{r}{4b^{2}})|\geq\left(\frac{r}{s}\right)^{\beta}
\frac{|Q(z_j,(2b^2+1)s)|}{C_7^m(8b^4+4b^2)^{\beta}}
\end{equation}
By (\ref{dbl}), $z$\ is in $Q(z_j,b^2s).$\ Now Lemma~\ref{psdmet}\  
implies: $Q(z_j,(2b^2+1)s)\supseteq
Q(z,(b^2+1)s).$\ This, (\ref{rro}) and (\ref{boi}) together imply:
\[
|Q(z,(b^2+1)s)|\,\geq\,p\left(\frac{r}{s}\right)^{\beta}
\frac{|Q(z,(b^2+1)s)|}{C_7^m(8b^4+4b^2)^{\beta}}\,.
\]
This proves (\ref{cota}) with $C_8=C_7^m(8b^4+4b^2)^{\beta}.$
\cqd

\begin{lemma}
\label{lc9}
There are constants $C_9,$\ and $\gamma$\ such that
\begin{equation}
\label{rs}
\frac{u(sQ)}{u(rQ)}\leq C_9\left(\frac{s}{r}\right)^{\gamma}
\ \mbox{and}\ \ 
\frac{v(sQ)}{v(rQ)}\leq C_9\left(\frac{s}{r}\right)^{\gamma}
\end{equation}
for every ``rectangle'' $Q$\ and for every $0<r<s.$
\end{lemma}

\prf
It follows from (\ref{dbl}) and (\ref{dob}) that, if $l$\ is an integer such 
that $b^2<2^l,$\ then $u(2Q)\leq C_3^{l+1}u(Q)$\ and $v(2Q)\leq C_4^{l+1}v(Q)$\ 
for all $Q.$\ Arguing similarly as for the proof of (\ref{invdob}), we can 
get (\ref{rs}) with $C_9=\max\{C_3^{l+1},C_4^{l+1}\},$\ 
and $\gamma=\log C_9/\log 2.$
\cqd

The following theorem plays here the role of Theorem~D in \cite{JF}. The 
explicit form of the constants in (\ref{lab}), valid for arbitrarily small
$\epsilon,$\ is needed for an efficient control of the constants that show up
in the iteration process.
%
%
\begin{thm}
\label{D}
Under the hypotheses of Theorem~\ref{main}, there are constants $\alpha$\ and 
$C_{10}$\ such that the estimate
\begin{equation}
\label{lab}
\frac{\epsilon^\alpha}{C_{10}}
\left[\frac{1}{u(Q)}\int_{Q}|g(z)|^q u(z)\dz\right]\leq 
\end{equation}
\[
\left[
\left(
\frac{s^2}{v(Q)}\int_{(1+\epsilon)Q}|\nl g(z)|^2v(z)\dz
\right)^{\frac{1}{2}}+
\left(
\frac{1}{u(Q)}
\int_{(1+\epsilon)Q} g(z)^2u(z)\dz
\right)^{\frac{1}{2}}
\right]^{q}
\]
holds for every $Q=Q(z,s),$\ for every $0<\epsilon<1,$\ and for every 
Lipschitz continuous function $g,$\  
where $q>2$\ is the constant provided by \ap{\bf SP}\fp.
\end{thm}

\prf
Let us apply Proposition~\ref{cover}\ with $r=\epsilon s/b^2$\ and let 
the $Q$'s then obtained be denoted by $Q_j=Q(z_j,r),$\ $j=1,\cdots,m.$ By 
(\ref{spc}) we get:
\[
\int_{Q}|g(z)|^q u(z)\dz\leq
\sum_{j=1}^{p}u(Q_j)
\left[ C_6r
\left(
\frac{1}{v(Q_{j})}\int_{b^{2}Q_{j}}|\nl g(z)|^2 v(z)\dz
\right)^{\frac{1}{2}}
\right.
\]
\begin{equation}
\label{tez}
+\,
\left.
\left(
\frac{1}{u(Q_{j})}\int_{b^2Q_j}g(z)^2u(z)\dz
\right)^{\frac{1}{2}}
\right]^q.
\end{equation}
By Lemma~\ref{psdmet}, we have $b^2Q_j\subseteq Q(z,s+b^2r),$\ 
and hence the integrals on $b^{2}Q_{j}$\ inside the brackets 
in (\ref{tez}) may be replaced
by integrals on $(1+\epsilon)Q.$\ We then estimate $u(Q)/u(Q_j)$\ and $v(Q)/v(Q_j)$\ using (\ref{rs}) and $Q(z,s)\subseteq Q(z_j,(b^2+1)s)$\ (which 
follows from Lemma~\ref{psdmet}\ and Proposition~\ref{doble}). This way we see 
that the expression between brackets in (\ref{tez}) is bounded by the expression between brackets in (\ref{lab}) times 
$C_9^{\frac{1}{2}}\max\{C_6,1\}[(b^2+1)s/r]^{\frac{\gamma}{2}}.$\ 
Next we use that $Q_j\subseteq 2Q$\ (which follows from Lemma~\ref{psdmet}), 
to get $u(Q_j)\leq C_9u(Q)$\ (by the proof of Lemma~\ref{lc9}).
After using (\ref{cota}), we finally get (\ref{lab}) with 
\[
C_{10}=C_8C_9^{\frac{2+q}{2}}\max\{C_6,1\}^q(b^4+b^2)^{\frac{q\gamma}{2}}
b^{2\beta}
\]
and $\alpha=\beta+q\gamma/2.$\
\cqd
%
%
%
%
%
%
%

\section{Moser iteration and Harnack inequality}\label{moserze}

We start this section with the construction of the test functions addapted
to rectangles mentioned in the Introduction. 

\begin{prop}
\label{test}
Given any $z_{\circ}\in\rN$\ and any $0<r_1<r_2,$\ there is a smooth function 
$\eta$\ equal to one everywhere on $Q(z_{\circ},r_1),$\ 
with support contained in 
$Q(z_{\circ},r_2),$\ and such that $0\leq\eta(z)\leq 1$\ and 
$|\nl\eta(z)|\leq C_{11}/(r_2-r_1)$\ for all $z\in\rN,$\ with $C_{11}$\ 
denoting the constant $2\sqrt{N}.$
\end{prop}

\prf
Choose $\psi$\ a smooth function on $\mathbb R$\ identical to one on 
$(-\infty,0],$\ with support contained in $(-\infty,1),$\
and such that $0\leq\psi(t)\leq 1$\
and $|\psi^{\prime}(t)|\leq 2$\ for all $t\in{\mathbb R}.$\ Given 
$z_{\circ}=(x_{\circ},y_{\circ})\in\rn\times\rrm$\ and $0<r_1<r_2,$\ define
\[
\eta(x,y)=\varphi\left(\frac{|x-x_{\circ}|}{r_2}\right)
          \varphi\left(\frac{|y-y_{\circ}|}{r_2\Lambda(z_{\circ},r_2)}\right),\ 
          (x,y)\in\rn\times\rrm\,,
\]
where $\varphi(t)=\psi(\frac{r_2t-r_1}{r_2-r_1}).$\ It is straightforward 
to check that this $\eta$\ does it. \cqd

\begin{defn} {\em An element $f\in\ho$\ is a weak subsolution of $Lf=0$\ if
$a_{\circ}(f,\theta)\leq 0$\ for all non-negative $\theta$\ in $\hoo.$}
\end{defn}

\begin{defn} {\em 
Given $M>0$\ and $d\geq 1,$\ let the function 
$\hmb$\ (continuously differentiable with bounded derivative)
be defined by $\hmb(t)=t^d$\ if $t\in[0,M],$\ and 
$\hmb(t)=M^d+d M^{d-1}(t-M)$\ if $t>M.$ }
\end{defn}

\begin{prop}
\label{groiss}
Let $f\in\hq,$\ $Q=Q(z_{\circ},h),$\ be a non-negative subsolution of 
$Lf=0$\ and let $f_k$\ be 
a sequence of non-negative Lipschitz continuous functions on 
$\overline{Q}$\ converging to $f$\ in $\hq.$\ 
Given $\frac{1}{2}\leq s<t\leq 1,$\
$M>0$\ and $\beta\geq 1,$\ there are a subsequence $\fkj$\ of $f_k$\ and 
a sequence  $\delta_j\geq 0$, $\delta_j\rightarrow 0,$\ such that 
for all $j$\ we have:
\begin{equation}
\label{jorge}
\int_{sQ}|\nl(\hmb\circ\fkj)|^2v\leq\delta_j+
\frac{4C_{11}^{2}}{(t-s)^{2}h^{2}}\int_{tQ}|\fkj\cdot(\hmb^{\prime}\circ\fkj)|^2u.
\end{equation}
\end{prop}

Proposition~\ref{groiss}\ can be given a proof
almost identical to the first part of the proof of Lemma (3.1) in \cite{ChW}\ 
(pages 1117 to 1119). One only needs to replace their Euclidean ball $B$\ by
our rectangle $Q,$\ and their ellipticity condition (1.1) by ours 
(\ref{ell}). When (\ref{ell}) is applied, our $\nl$\ will show up, replacing
their $\nabla.$\ Also, one should take 
$\eta$\ as the test function constructed in Proposition~\ref{test},\ with 
$r_1=hs$\ and $r_2=ht.$\ Since the support of the chosen $\eta$\ is contained
in the open set $Q(z_{\circ},ht),$\ we may allow $t$\ to be  
equal to one (this fact is needed in our iteration).

An inequality to be derived from (\ref{lab}) and (\ref{jorge}) will be 
iterated in the proof of the next proposition, which corresponds to a 
weaker version of Lemma (3.1) in \cite{ChW}. 

\begin{prop}
\label{crucial}
If $f\in\hq,$\ $Q=Q(z_{\circ},h),$\ is a non-negative subsolution of 
$Lf=0,$\ then the estimate
\begin{equation}
\label{crcl}
\left(
\mbox{\tt ess\,sup}_{aQ}\tilde{f}
\right)^p\,\leq\,
\frac{C_{12}}{(1-a)^\delta}[p\mu(Q)]^{\frac{2q}{q-2}}\frac{1}{u(Q)}
\int_{Q}\tilde{f}^pu
\end{equation}
holds for every $a\in[\frac{1}{2},1)$\ and every $p\geq 2,$\ with 
$\delta$\ and $C_{12}$\ denoting constants explicitly defined below
\ap at the end of the proof\,\fp, and
$\mu(Q)=u(Q)^{\frac{1}{2}}v(Q)^{-\frac{1}{2}}.$\ \ap We recall that $q$\ 
arises in \ap{\bf SP}\fp .\fp
\end{prop}

\prf
Given $\frac{1}{2}\leq s<t\leq1$\ and $d\geq 1,$\ let us first use
Proposition~\ref{groiss}\ to extract a subsequence $\fkj$\ of 
a sequence $f_k$\ of non-negative Lipschitz continuous functions on 
$\overline{Q}$\ converging to $f$\ in $\hq$\ for which (\ref{jorge})
is true. Then, let us apply (\ref{lab}) 
to the rectangle $sQ,$\ for some $\epsilon$\ satisfying $(1+\epsilon)s<t$\ 
and for $g=\hmb\circ\fkj.$\ Then let us apply (\ref{jorge}) with 
$(1+\epsilon)s$\ replacing $s.$\ Next, we use (\ref{rs}) with $t$\ and $s$\ 
replacing $s$\ and $r,$\ respectively, and taking advantage of the fact that 
$1<t/s\leq 2.$\  Finally, after using that 
$\hmb(\varphi)\leq \varphi\hmb^\prime(\varphi)$\ 
for all $\varphi\in{\mathbb R},$\ we get
\begin{equation}
\label{mess}
\left[
\frac{1}{u(sQ)}\int_{sQ}|\hmb\circ\fkj|^qu
\right]^{\frac{1}{q}}
\leq
\frac{C_{10}^{\frac{1}{q}}\epsilon^{-\frac{\alpha}{q}}sh}
{v(sQ)^{\frac{1}{2}}}\delta_j^{\frac{1}{2}}\,+\,
2^{\frac{\gamma}{2}}C_{9}^{\frac{1}{2}}C_{10}^{\frac{1}{q}}
\epsilon^{-\frac{\alpha}{q}}\cdot
\end{equation}
\[
\left[
2C_{11}\mu(sQ)\frac{s}{t-(1+\epsilon)s}+1
\right]
\left[
\frac{1}{u(tQ)}\int_{tQ}|\fkj\cdot(\hmb^\prime\circ\fkj)|^2u
\right]^{\frac{1}{2}}.
\]

Now we want to let $j$ first, and then $M,$\ go to infinity. We may suppose,
passing to another subsequence if necessary, that $\fkj$\ converges to 
$\tilde{f}$\ pointwise, almost everywhere with respect to the measure 
$u(z)\dz.$\ Using Fatou's Lemma on the left-hand side and Lebesgue's 
convergence theorem on the right (again, this is the same argument as 
Chanillo and Wheeden's, on page 1120 of \cite{ChW}), one can see that it 
is legitimate to replace $\fkj$\ by $\tilde{f}$\ in (\ref{mess}), and then
$\hmb\circ\tilde{f}$\ by $\tilde{f}^{d}$\ and
$\hmb^{\prime}\circ\tilde{f}$\ by $\tilde{f}^{d-1}.$

Since $\frac{1}{2}\leq s<(1+\epsilon)s<t\leq 1,$\ then $s/[t-(1+\epsilon)s]$\ is greater than one. By (\ref{ell}), it follows that $\mu(sQ)\geq 1.$\ Hence, 
the ``+1'' inside the first pair of brackets at the right-hand side of 
the inequality (\ref{mess}) may 
be absorbed by the constant at its left, which will then be multiplied by 
two. Next we raise to the $\frac{1}{d}$-th power both sides of the 
inequality and change notation, writing $r=2d$\ and $q=2\sigma.$\ After all 
that is taken into account, we will have deduced from (\ref{mess}) the 
estimate
\begin{equation}
\label{moser}
\left[
\frac{1}{u(sQ)}\int_{sQ}\tilde{f}^{r\sigma}u
\right]^{\frac{1}{r\sigma}}\leq
\left[
\frac{C_{13}\epsilon^{-A}\mu(sQ)rs}{t-(1+\epsilon)s}
\right]^{\frac{2}{r}}
\left[
\frac{1}{u(tQ)}\int_{tQ}\tilde{f}^{r}u
\right]^{\frac{1}{r}},
\end{equation}
for all $r\geq 2,$\ 
with $C_{13}=2^{\frac{2+\gamma}{2}}C_{9}^{\frac{1}{2}}C_{10}^{\frac{1}{q}}
C_{11}$\ and $A=\frac{\alpha}{q}.$\ 

Let $a\in[\frac{1}{2},1)$\ and $p\geq 2$\ be given and define 
$a_j=a+(1-a)/(j+1).$\ For each non-negative integer $j,$\ let us apply 
(\ref{moser}) with $t=a_j,$\ $s=a_{j+1},$\ 
$\epsilon=\epsilon_j=(a_{j+1}-a_{j+2})/a_{j+1}$\ and $r=\sigma^jp.$\
Let us apply (\ref{moser}) again to the right-hand side of the
inequality thus obtained, but with $t=a_{j-1},$\ $s=a_j,$\ 
$\epsilon=\epsilon_{j-1}$\ and $r=\sigma^{j-1}p.$\ 
By repeating this procedure, after $j+1$\ steps we will get:
\begin{equation}
\label{moserj}
\left[
\frac{1}{u(a_{j+1}Q)}\int_{a_{j+1}Q}\tilde{f}^{p\sigma^{j+1}}u
\right]^{\frac{1}{p\sigma^{j+1}}}
\leq
\end{equation}
\[
\left\{
\prod_{k=0}^{j}
\left[
\frac{C_{13}\epsilon_k^{-A}\mu(a_{k+1}Q)p
\sigma^ka_{k+1}}{a_{k}-(1+\epsilon_{k})a_{k+1}}
\right]^{\frac{2}{p\sigma^{k}}}
\right\}
\left[\frac{1}{u(Q)}\int_{Q}\tilde{f}^pu
\right]^{\frac{1}{p}}.
\]
Since $a<a_{j+1}<2a$\ for all $j,$\ it follows from Lemma~\ref{lc9}\
that the left-hand side of (\ref{moserj}) is greater than or equal to
\[
\left[
\frac{2^{\gamma}C_{9}}{u(aQ)}\int_{aQ}\tilde{f}^{p\sigma^{j+1}}u
\right]^{\frac{1}{p\sigma^{j+1}}},
\]
which converges to $\mbox{\tt ess\,sup}_{aQ}\tilde{f}$\ as $j$\ tends to 
infinity. On the right-hand side of (\ref{moserj}) we may replace 
$\mu(a_{j+1}Q)$\ by $\sqrt{2^{\gamma}C_9}\mu(Q),$\ due to 
Lemma~\ref{lc9}. Hence, all we need is to find a 
precise estimate for the product
\begin{equation}
\label{infpro}
\prod_{k=0}^{\infty}
\left[ C_{13}\sqrt{2^{\gamma}C_9}\mu(Q)p\sigma^{k}
\right]^{\frac{2}{p\sigma^{k}}}
\prod_{k=0}^{\infty}
\left[
\frac{\epsilon_k^{-A}a_{k+1}}{a_{k}-(1+\epsilon_{k})a_{k+1}}
\right]^{\frac{2}{p\sigma^{k}}}.
\end{equation}
The first of these products equals 
$[ C_{13}\sqrt{2^{\gamma}C_{9}}\mu(Q)p
\sigma^{\frac{1}{\sigma -1}}
]^{\frac{2\sigma}{p(\sigma -1)}}.$\ 
The second expression between brackets is equal to the left side of:
\[
\frac{(k+1)(k+3)^{A+1}(ak+a+1)^{A+1}}{2(1-a)^{A+1}}\leq
\frac{(k+3)^{2A+3}}{2(1-a)^{A+1}}.
\]
Hence, the second product in (\ref{infpro}) is bounded by 
\[
2^{-\frac{2\sigma}{p(\sigma-1)}}(1-a)^{-\frac{2(A+1)\sigma}{p(\sigma-1)}}
\left[\exp\sum_{k=0}^{\infty}\frac{\log(k+3)}{\sigma^{k}}\right]^{\frac{4A+6}{p}}
.\]
Defining
$
\delta=\frac{2(A+1)\sigma}{\sigma-1}
$\ and
\[ 
C_{12}=
\left[
\sqrt{2^{\gamma-2}C_{9}}C_{13}\sigma^{\frac{1}{\sigma-1}}
\right]^{\frac{2\sigma}{\sigma-1}}
\left[\exp
\sum_{k=0}^{\infty}\frac{\log(k+3)}{\sigma^k}
\right]^{4A+6}
\]
finishes the proof. \cqd

\begin{prop}
\label{groiss2}
Let $f\in\hq,$\ $Q=Q(z_{\circ},h),$\ be a strictly positive \ap
$f\geq\epsilon_{\circ}>0$\fp\ solution of $Lf=0$\ and let $f_k,$\
$f_k\geq\epsilon_{\circ},$\ be a sequence in $\lip(\overline{Q})$\ converging to $f$\ in $\hq.$\ 
Given $\frac{1}{2}\leq s<t\leq 1$\ and $\beta\leq 1,$\ with 
$-1\neq\beta\neq 0,$\ there are a subsequence $\fkj$\ of $f_k$\ and 
a sequence of non-negative reals $\delta_j\rightarrow 0$\ such that 
for all $j$\ we have:
\begin{equation}
\label{jorge2}
\int_{sQ}|\nl(\fkj^{\frac{\beta+1}{2}})|^2v\leq\delta_j+
\frac{(\beta+1)^2}{\beta^2}\frac{C_{11}^{2}}{(t-s)^{2}h^{2}}
\int_{tQ}\fkj^{\beta+1}u .
\end{equation}
\end{prop}

A proof for Proposition \ref{groiss2}\ can be given following exactly the 
same steps as in the first half of the proof of Lemma (3.11) in \cite{ChW}, 
pages 1121 and 1122, making the addaptations already described after the 
statement of Proposition~\ref{groiss}.

The proof of the following proposition follows the steps of 
Lemma (3.11) of \cite{ChW}, for $p<0$\ or $p\geq 2.$\ 
For $0<p<2,$\ we use a technique of Hardy and Littlewood, 
as in Lemma (3.17) of \cite{GW}.

\begin{prop}
\label{crucial2}
If $f\in\hq,$\ $Q=Q(z_{\circ},h),$\ is a non-negative solution of 
$Lf=0,$\ then the estimate
\begin{equation}
\label{crcl2}
\left(
\mbox{\tt ess\,sup}_{aQ}\tilde{f}
\right)^p\,\leq\,
\frac{C_{14}}{(1-a)^{\delta}}[1+|p|\mu(Q)]^{\frac{2q}{q-2}}
\frac{1}{u(Q)}\int_{Q}\tilde{f}^pu
\end{equation}
holds for every $a\in[\frac{1}{2},1)$\ and every $0\neq p\in{\mathbb R},$\ 
with $\delta$\ and $\mu(Q)$\ as defined in Proposition~\ref{crucial}\ and 
$C_{14}$\ denoting the constant explicitly defined below, at the end of 
the proof.
\end{prop}

\prf We may suppose that $f\geq\epsilon_{\circ}>0$\ and later let 
$\epsilon_{\circ}$\ tend to zero, as long as we make sure that none of 
the constants depends on $\epsilon_{\circ}.$\ 

Given $\beta\leq 1,$\ $-1\neq\beta\neq 0,$\ $\frac{1}{2}\leq s<t\leq 1$\
and $\epsilon >0$\ such that $(1+\epsilon)s<t,$\ we may combine (\ref{lab})
for $g=\fkj^{\frac{\beta+1}{2}}$\ and (\ref{jorge2}), and then let 
$j$\ go to infinity. 
Similarly as just before (\ref{moser}), with $r=\beta+1$\ and $\sigma=q/2,$\ we get:

\begin{equation}
\label{moser2}
\left[
\frac{1}{u(sQ)}\int_{sQ}\tilde{f}^{r\sigma}u
\right]^{\frac{1}{|r|\sigma}}\leq
\end{equation}
\[
\left(
\frac{C_{13}\epsilon^{-A}}{2}
\right)^{\frac{2}{|r|}}
\left[
\frac{|r|}{|r-1|}\cdot
\frac{s\mu(sQ)}{t-(1+\epsilon)s}+1
\right]^{\frac{2}{|r|}}
\left[
\frac{1}{u(tQ)}\int_{tQ}\tilde{f}^{r}u
\right]^{\frac{1}{|r|}}
\]
for all $r\leq 2,$\ $0\neq r\neq 1.$\ 

Now let $a\in[\frac{1}{2},1)$\ and $p<0$\ be given and let $a_j$\
and $\epsilon_j$\ be defined as in the proof of Proposition~\ref{crucial}. 
For each integer $j,$\ let us then apply (\ref{moser2}) with $r=\sigma^kp,$\ 
$t=a_k,$\ $s=a_{k+1}$\ and $\epsilon=\epsilon_k,$\ for $k=0,1,\cdots,j.$\ 
Iterating the $j+1$\ inequalities just obtained and letting $j$\
tend to infinity, similarly as before, we get:
\begin{equation}
\label{infpro2}
\mbox{\tt ess\,sup}_{aQ}\,\tilde{f}^{-1}
\leq K_{0}
\left[\frac{1}{u(Q)}\int_{Q}\tilde{f}^pu
\right]^{\frac{1}{|p|}},
\end{equation}
with
\[
K_{0}=\prod_{k=0}^{\infty}
\left[
\frac{C_{13}\epsilon_k^{-A}}{2}
\left(
\frac{\sigma^k|p|\mu(a_{k+1}Q)a_{k+1}}
{|\sigma^kp-1|\cdot[a_{k}-(1+\epsilon_{k})a_{k+1}]}
+1
\right)
\right]^{\frac{2}{|p|\sigma^{k}}}.
\]
Since at this point we are assuming $p<0,$\ we have $|\sigma^kp-1|\geq 1$\ 
for all $k\geq 0.$\ Taking into account also that $1\leq\sigma^k,$\ 
that $1\leq a_{k+1}/[a_k-(1+\epsilon_k)a_{k+1}],$\ that 
$\mu(a_{k+1}Q)\leq\sqrt{2^{\gamma}C_9}\mu(Q)$\ 
and that $1\leq\sqrt{2^{\gamma}C_9},$\ the infinite product above is seen 
to be bounded by:
\[
K_{0}\leq
\prod_{k=0}^{\infty}
\left[
\frac{C_{13}\sqrt{2^{\gamma}C_9}}{2}(1+|p|\mu(Q))\sigma^k
\frac{\epsilon_k^{-A}a_{k+1}}{a_{k}-(1+\epsilon_{k})a_{k+1}}
\right]^{\frac{2}{|p|\sigma^{k}}}.
\]
We may here use the estimates obtained at the end of the proof of 
Proposition~\ref{crucial} to conclude that (\ref{crcl2}) holds, if there we 
replace $C_{14}$\ by $C_{12}.$\ It follows from Proposition~\ref{crucial}\ 
that the same is true for $p\geq 2.$\ 

In case $0<p<2,$\ we have
$\sigma^jp$\ tending to infinity, but smaller than two for some values 
of $j.$\ Let us first suppose that $\sigma^{k}p\neq 1,$\
for every integer $k\geq 0.$\ Let then $l$\ be the integer such that
$\sigma^{l}p<2\leq\sigma^{l+1}p.$\ We may iterate as before, but using 
(\ref{moser2}) at the first $l+1$\ steps of the iteration and (\ref{moser}) 
after that. We get: 
\[
\mbox{\tt ess\,sup}_{aQ}\tilde{f}
\,\leq\,K_{1}
\left[
\frac{1}{u(Q)}\int_{Q}\tilde{f}^pu
\right]^{\frac{1}{p}},
\]
with 
\[
K_{1}=\prod_{k=0}^{l}
\left[
\frac{C_{13}\epsilon_k^{-A}}{2}
\left(
\frac{\sigma^kp\mu(a_{k+1}Q)a_{k+1}}
{|\sigma^kp-1|\cdot[a_{k}-(1+\epsilon_{k})a_{k+1}]}
+1
\right)
\right]^{\frac{2}{p\sigma^{k}}}\cdot
\]\[
\prod_{k=l+1}^{\infty}
\left[
\frac{C_{13}\epsilon_k^{-A}\sigma^kp\mu(a_{k+1}Q)
a_{k+1}}{a_{k}-(1+\epsilon_{k})a_{k+1}}
\right]^{\frac{2}{p\sigma^{k}}}.
\]

In order to get a good estimate for $K_{1},$\ let us further suppose 
that $p=\sigma^{j}(\sigma+1)/2,$\ for some $j\in{\mathbb Z}.$\ Then it will 
hold that $|\sigma^{k}p-1|\geq(\sigma-1)/(2\sigma),$\ for every integer
$k\geq 0.$\ We may proceed as we did for the other infinite products, using
in adition that $1<2\sigma/(\sigma-1),$\ and prove that 
(\ref{crcl2}) holds for these values of $p,$\ with $C_{14}$\ replaced by 
$C_{15}=C_{12}(\frac{2\sigma}{\sigma-1})^{\frac{2\sigma}{\sigma-1}}.$\  

By Remark~\ref{rstr}, we may apply the result we have just obtained 
with $\alpha Q$\ replacing $Q,$\ for any $\alpha\in(0,1).$\ Given
$\frac{1}{2}\leq{\alpha^\prime}< \alpha\leq 1$\ and $\pb$\ 
belonging $X=\{\sigma^{j}(\sigma+1)/2;\,j\in{\mathbb Z}\},$\ we get:
\begin{equation}
\label{putz}
\left(
\mbox{\tt ess\,sup}_{\alpha^{\prime}Q}\tilde{f}
\right)^{\pb}\,\leq\,
\frac{C_{15}(2^{\gamma}C_{9})^{\frac{q}{q-2}}}
{(\alpha-\alpha^{\prime})^{\delta}}
[1+\pb\mu(Q)]^{\frac{2q}{q-2}}
\frac{1}{u(\alpha Q)}\int_{\alpha Q}\tilde{f}^{\pb}u,
\end{equation}
where we have used Lemma~\ref{lc9}\ and $1\leq C_{9}$\ in order to  
replace $\mu(\alpha Q)$\ by $\mu(Q)$\ inside the brackets. 

Let us define 
\[
I_p=\frac{[1+p\mu(Q)]^{\frac{2q}{2q-2}}}{u(Q)}\int_{Q}\tilde{f}^pu, 
\ \ p\in(0,2),
\]
and $E(\alpha)=\mbox{\tt ess\,sup}_{\alpha Q}\tilde{f}.$\ 
Given any $p\in (0,2)\setminus X,$\ let $\pb\in X$\ be such that 
$\frac{\pb}{\sigma}<p<\pb.$\ By Lemma~\ref{lc9}\ and (\ref{putz}), we get
\begin{equation}
\label{tqsdi}
E(\alpha^{\prime})^{\pb}\leq 
\frac{C_{16}}{(\alpha-\alpha^{\prime})^{\delta}}
E(\alpha)^{\pb-p}I_p,
\end{equation}
with $C_{16}=C_{15}(2^{\gamma}C_{9})^{\frac{2q-2}{q-2}}
\sigma^{\frac{2q}{q-2}}.$

Given $a\in[\frac{1}{2},1),$\ let $\alpha_{k}$\ be a strictly increasing 
sequence such that $\alpha_{0}=a,$\ and $\lim \alpha_k< 1.$\ Let us take
the logarithm of (\ref{tqsdi}) and iterate, with 
$\alpha^{\prime}=\alpha_k$\ and $\alpha=\alpha_{k+1},$\ $k=0,1,\cdots$.
With $\theta=(\pb-p)/\pb,$\ we get
\[
\log E(a)\leq\frac{1}{\pb}\sum_{k=0}^{\infty}\theta^k
\log\frac{C_{16}}{(\alpha_{k+1}-\alpha_{k})^{\delta}}
\]
\begin{equation}
\label{log}
+\limsup_{k\rightarrow\infty}\theta^{k+1}\log E(\alpha_{k+1})
+\frac{1}{\pb(1-\theta)}\log I_p\,;
\end{equation}
noting that, since $C_{16}>1$\ and $\alpha_{k+1}-\alpha_{k}<1/2,$\
the terms of the series in the above inequality are positive.

It follows from Proposition~\ref{crucial}\ for $p=2$\ that 
$E(\lim\alpha_k)$\ is finite. Since $\theta<1,$\ we then get 
$\limsup_{k\rightarrow\infty}\theta^{k+1}\log E(\alpha_{k+1})=0.$\ 
To estimate the sum in (\ref{log}), we need to make a 
precise choice of $\alpha_k.$\ If we let
\[
\alpha_{k}=a+(1-a)\frac{\sum_{j=1}^{k}j^{-2}}{2\sum_{j=1}^{\infty} j^{-2}},\ 
k\geq 1,
\]
we get $\alpha_{k+1}-\alpha_k\geq (1-a)/[4(k+1)^2].$\ 
Since $\pb(1-\theta)=p,$\ we get:
\[
p\log E(a)\leq\log\frac{C_{16}4^{\delta}}{(1-a)^{\delta}}
+\sum_{k=0}^{\infty}\theta^{k}\log(k+1)^{2\delta}+\log I_p.
\]
Exponentiating both sides of this inequality and defining
\[
C_{14}=\max\left\{C_{12},\,C_{15},\,
4^{\delta}C_{16}\exp\left[\sum_{k=0}^{\infty}\theta^{k}\log(k+1)^{2\delta}
\right]\right\}
\]
finishes the proof. \cqd

\begin{prop}\label{p1}
Let $\Omega$\ be a bounded open set of $\rN,$\ and let $f\in\ho$\ be a 
positive weak solution of $Lf=0,$\ bounded below by a positive number.
Let $z_{0}\in\Omega$\ and $h>0$\ be such that $bB\subseteq\Omega,$\ where 
$B=B(z_{0},h).$\ For each $\alpha\in [\frac{1}{2},1),$\ define 
$k(\alpha,f)$\ by $\log k(\alpha,f)=m_{\alpha bB}(\log\tilde{f}).$\ \ap See 
page~{\em \pageref{pagina}}\fp
Then there is a constant $C_{17}$\ such that, if $z_{0}$\ and $h$\ are 
such that $b^2Q\subseteq\Omega,$\ where $Q=Q(z_{0},h),$\ then the inequality
\begin{equation}
\label{pnulo}
u(\{x\in\alpha Q;\,|\log\frac{\tilde{f}(x)}{k(\alpha,f)}|>t\})
\leq\frac{C_{17}\mu(Q)u(\alpha Q)}{(1-\alpha)t}
\end{equation}
holds for every $t>0$\ and every $\alpha\in[\frac{1}{2},1).$
\end{prop}

\prf
This proposition can be given a proof very similar to that of Lemma (3.13)
of \cite{ChW}. We are going to highlight a few points, refering
to Chanillo and Wheeden's article for more details.

Let $f_k$\ denote a sequence of positive Lipschitz continuous functions, 
uniformly bounded away from zero, converging to $f$\ in $\ho.$\ With the 
aid of the test function $\eta$\ (built in Proposition~\ref{test}\ 
-- here we take $r_{1}=\alpha h$\ and $r_{2}=h$), we can extract 
from $f_k$\ a subsequence, which we will still denote by $f_k,$\ such that
\begin{equation}
\label{vtr}
\int_{\alpha Q}|\nl (\log f_k)|^2v\leq
\frac{4C_{11}^{2}u(Q)}{(1-\alpha)^{2}h^{2}}+\delta_{k},
\end{equation}
for some $\delta_k\rightarrow 0.$

With $g=\log f_k,$\ let us apply (\ref{spc0}) with $q$\ replaced by 2 (this 
is allowed by H\"older's inequality) and $Q$\ replaced by $\alpha Q$. Next, 
let us apply (\ref{vtr}) with $Q$\ replaced by $b^2Q.$ Using also 
Lemma~\ref{lc9}, we get:
\begin{equation}
\label{dica}
\int_{\alpha Q}|\log(f_k)-m_{\alpha bB}(\log f_k)|^2u
\leq
\frac{C_{17}^{2}}{(1-\alpha)^{2}}\mu(Q)^2u(\alpha Q)+\delta^{\prime}_{k},
\end{equation}
with $C_{17}=2C_{6}C_{9}^{\frac{1}{2}}C_{11}b^{\gamma-2}$\ and 
$\delta_{k}^{\prime}\rightarrow 0.$\ Using that $f_k$\ is uniformly bounded
away from zero, one can see that the $\lim_k$\ of the left-hand side of 
(\ref{dica}) is equal to 
$\int_{\alpha Q}|\log\tilde{f}-\log k(\alpha,f)|^2u.$\ The Proposition
now follows from Chebyshev's and Cauchy-Schwartz's inequalities.\cqd

The following lemma for $w\equiv 1$\ is essentially Lemma~3 of \cite{M1},\ 
whose proof also works for the case of an arbitrary weight $w.$\

\begin{lemma}
\label{bombieri}
\ap {\bf Bombieri-Moser}\fp\ 
Let $w$\ be a \ap non-negative\fp\ weight on $\rN,$\ and let $f$\ be a 
bounded non-negative measurable function defined on a bounded measurable set 
$E$. Suppose there is a family $E_t,$\ $t\in(0,1],$\ of measurable sets with 
$w(E_{t})>0$\ for all $t,$\ $E_1=E$\ and $E_s\subset E_t$\ if $s<t.$\ 
Assume there are $\mu,c,d>0,$\ such that 
\begin{equation}
\label{b1}
\mbox{\tt ess\,sup}_{E_{s}}\,f^p\leq\frac{c}{(t-s)^{d}}
\frac{1}{w(E_{1})}\int_{E_{t}}f^pw,
\end{equation}
for all $p,$\ $s,$\ and $t$\ such that $0<p<\mu^{-1}$\ and 
$\frac{1}{2}\leq s<t\leq 1;$\ and
\begin{equation}
\label{b2}
w(\{x\in E_{1};\,\log f(x)>\tau \} )\leq\frac{c\mu}{\tau}w(E_{1}),
\end{equation}
for all $\tau >0.$\ Then there exists
$C>0,$\ depending only on $c$, such that
\begin{equation}\label{35}
\mbox{\tt ess\,sup}_{E_{\alpha}}\,f
\leq \exp\left[\frac{C\mu}{(1-\alpha)^{2d}}
\right]\,,
\end{equation}
for all $\alpha\in [\frac{1}{2},1).$
\end{lemma}

\noindent{\bf Proof of Theorem~\ref{main}:}\ We may suppose that $\tilde{f}$\
is bounded away from zero, otherwise we could add an $\epsilon > 0$\ and
later let $\epsilon\rightarrow 0.$\ 

Let $B=B(z_{0},h)$\ be such that $2b^3B\subseteq\Omega$\ and let 
$Q=Q(z_{0},h).$\ With $w=u$\ and $E_{t}=\frac{3t}{2}Q,$\ we are going to 
apply Lemma~\ref{bombieri}\ to the functions 
$\tilde{f}/k$\ and $k/\tilde{f},$\ where 
$k=\exp[m_{\frac{3}{2}bB}(\log\tilde{f})]$. Notice that
$\tilde{f}$\ is bounded on $E_1,$\ since the closure of $E_1$\ is contained
in $\Omega,$\ and we may then apply Proposition~\ref{crucial}\ with $p=2$\ 
for a rectangle slightly larger than $E_1 .$\ Choosing, for example, 
\[
c=\max\{4C_{17}\sqrt{2^{\gamma}C_{9}},2^{\gamma}C_{9}C_{14}
(2^{\frac{\gamma}{2}+1}C_{9})^{\frac{2q}{q-2}}\},
\]
we can check that (\ref{b1}) and (\ref{b2}) with $d=\delta$\ and 
$\mu=\mu(Q)$\ hold for both
$\tilde{f}/k$\ and $k/\tilde{f},$\ by Proposition~\ref{crucial2}\ and 
Proposition~\ref{p1}, and by also using that $u$\ is 
doubling (Lemma~\ref{lc9}). We remark that $2b^{3}B\subseteq\Omega$\
implies that $2b^{2}Q\subseteq\Omega,$\ and we may apply (\ref{pnulo})
with $2Q$\ replacing $Q$\ and $\alpha=\frac{3}{4}.$\ Choosing 
$\alpha=\frac{2}{3}$\ in (\ref{35}) for $\tilde{f}/k$\ and $k/\tilde{f},$\ 
we see that $\mbox{\tt ess\,sup}_{Q}(\tilde{f}/k)$\  and  
$\mbox{\tt ess\,sup}_{Q}(k/\tilde{f})=
[\mbox{\tt ess\,inf}_{Q}(\tilde{f}/k)]^{-1}$\ are both bounded by 
$\exp(3^{2d} C\mu).$\ Taking the product of these two inequalities, we get:
\begin{equation}
\label{36}
\mbox{\tt ess\,sup}_{Q}\tilde{f}\,\leq\,\exp(2C3^{2d}\mu)
\mbox{\tt ess\,inf}_{Q}\tilde{f}.
\end{equation}

Now let $B=B(z_{0},h)$\ be such that $2b^4\subseteq\Omega.$\ We may apply 
(\ref{36}) for the rectangle $bQ.$\ By (\ref{dbl}) we thus have
\[
\mbox{\tt ess\,sup}_{B}\tilde{f}\,\leq\,
\mbox{\tt ess\,sup}_{bQ}\tilde{f}\,\leq\,          \]\[
\exp[2C3^{2d}\mu(bQ)]\mbox{\tt ess\,inf}_{bQ}\tilde{f}
\,\leq\,\exp[2C3^{2d}\mu(bQ)]\mbox{\tt ess\,inf}_{B}\tilde{f}.
\]
This proves (\ref{harn}) with $K=2C3^{2d}$\ but with $\mu=\mu(bQ)$\
instead of $\mu=u(B)^{\frac{1}{2}}v(B)^{-\frac{1}{2}}.$\ Since $u$\
and $v$\ are doubling, those two quantities are comparable.\cqd

\section*{Acknowledgements}
We thank Richard Wheeden for suggesting the problem. The three authors were 
partially supported by Brazilian agency CNPq.

\vskip0.4cm

2000 {\em Mathematics Subject Classification}: 35J70.

\vskip0.4cm

{\footnotesize{\sc

Instituto de Matem\'atica e Estat\'\i stica

Universidade de S\~ao Paulo

Caixa Postal 66281, 05311-970 S\~ao Paulo SP, Brazil.}

\texttt{diniz@ime.usp.br}

\vskip0.3cm {\sc

Facultad de Ingenieria

Universidad de la Rep\'ublica 

Julio Herrera y Reissig 567, Montevid\'eo, C.P. 11300, Uruguay.}

\texttt{jorgeg@fing.edu.uy}

\vskip0.3cm {\sc

Instituto de Matem\'atica e Estat\'\i stica

Universidade de S\~ao Paulo

Caixa Postal 66281, 05311-970 S\~ao Paulo SP, Brazil.}

\texttt{melo@ime.usp.br}}

\end{document}